\documentclass[12pt,thmsa]{article}
\usepackage{amsmath}
\usepackage{graphicx}
\usepackage{amsfonts}
\usepackage{amssymb}
\usepackage{mathrsfs}

\newcommand{\X}{\mathcal X}
\newcommand{\Y}{\mathcal Y}

\newcommand{\HX}{\mathcal H_{\X}}
\newcommand{\HY}{\mathcal H_{\Y}}

\newtheorem{prop}{Proposition}
\newtheorem{lem}{Lemma}
\newtheorem{thm}{Theorem}
\newtheorem{corl}{Corollary}

\begin{document}

\begin{center}
{\Large \textbf{Testing independence of functional variables by an Hilbert-Schmidt independence criterion estimator}}

\bigskip

Terence Kevin  MANFOUMBI DJONGUET, Alban MBINA MBINA and  Guy Martial  NKIET 

\bigskip

\textsuperscript{}URMI, Universit\'{e} des Sciences et Techniques de Masuku,  Franceville, Gabon.

\bigskip

E-mail adresses : tkmpro95@gmail.com,   alban.mbinambina@univ-masuku.org, guymartial.nkiet@univ-masuku.org.

\bigskip
\end{center}

\noindent\textbf{Abstract.} We propose an estimator of   the Hilbert-Schmidt Independence Criterion obtained from an appropriate  modification of the usual estimator. We then get asymptotic normality of this estimator both under independence hypothesis and under the alternative hypothesis. A  new  test for independence  of  random variables valued into metric spaces is then  introduced, and a   simulation study that allows to compare   the proposed test   to an  existing one is provided.

\bigskip

\noindent\textbf{AMS 1991 subject classifications: }62E20, 46E22.

\noindent\textbf{Key words:} Asymptotic normality;  Kernel method;   Hilbert-Schmidt independence criterion; Reproducing kernel Hilbert space;  Functional data analysis.
\section{Introduction}

\noindent Functional data analysis is a branch of statistics allowing   to set up methods for processing data in the form of curves representing, for example, the evolution of random phenomena over time. It have known, over the last two decades, a  very important development which  allowed the introduction of methods for studying  links between functional variables through functional regression models, non-correlation tests (e.g., Kokoszka et al. (2008), Aghoukeng Jiofack and Nkiet (2010)), independence tests (e.g., G\'orecki et al. (2020), Lai et al. (2021)). However, only a few works focus on the problem of independence testing between functional variables. Recently, Lai et al. (2021)  introduced the angle covariance to characterize independence and proposed  a permutation based test whereas  G\'orecki et al. (2020)  proposed the use of the  Hilbert-Schmidt Independence Criterion (HSIC) for testing independence  from   multivariate functional data. HSIC, introduced by Gretton et  al. (2005),   is one of the most successful non-parametric dependence measure defined for random variables with values  in metric spaces, allowing to measure the dependency between univariate or multivariate random variables, but also random variables valued into  more complex structures such as structured, high-dimensional or functional   data. For testing independence, an empirical   estimator of HSIC was proposed  in Gretton (2005) as test statistic; but its  asymptotic distribution under null hypothesis    is an infinite sum of distributions   (see Zhang et al. (2018)) and, consequently,  can not been used for performing the test,  so requiring the use of methods such as permutation method to calculate  $p$-value. For overcoming   this  drawback, we  adopt  in this paper  an approach introduced in  Makigusa and Naito (2020), and used in Balogoun et al. (2021),    consisting in constructing a test statistic  from   an  appropriate modification of  the aforementioned estimator   in order to yield asymptotic normality both under the  null hypothesis and under the alternative. This approach   allows us to propose a new independence test for random variables  valued into  metric spaces, including functional variables. The rest of the paper is organized as follows.     The HSIC  is recalled in Section 2, and Section 3 is devoted to its estimation by a modification of the naive estimator,  and to the main results.  A simulation study on functional data that allows to compare   the proposed test   to that of Lai et al. (2021) is  given in Section 4. All the proofs are postponed in Section 5.

\section{HSIC and independence test}

\noindent Let $X$ and $Y$ be two random variables defined on a probability space $(\Omega,\mathscr{A},\mathbb{P})$ and with values in metric spaces $\mathcal{X}$ and  $\mathcal{Y}$ respectively. We consider reproducing kernel Hilbert spaces   $\HX$ and  $\HY$  of functions  from  $\X$ and $\Y$, respectively, to $\mathbb{R}$  with associated   kernels    $K\,:\,\mathcal{X}^2\rightarrow\mathbb{R}$ and $L\,:\,\mathcal{Y}^2\rightarrow\mathbb{R}$. They are symmetric functions such that, for any $(f,g)\in \HX\times\HY$ and any $(x,y)\in\mathcal{X}\times\mathcal{Y}$, one has $K(x,\cdot)\in\HX$, $L(y,\cdot)\in\HY$, $f(x)=\langle K(x,\cdot),f\rangle_{\HX}$ and  $g(y)=\langle L(y,\cdot),g\rangle_{\HY}$ (see Berlinet and Thomas-Agnan, 2004), where  $\langle \cdot , \cdot \rangle_{\HX}$  (resp. $\langle \cdot , \cdot \rangle_{\HY}$) denotes the inner product of $\HX$ (resp. $\HY$).   Throughout this paper, we assume that $K$ and $L$ satisfy the following assumption:

\bigskip

\noindent $(\mathscr{A}_1):$  $\Vert  K\Vert_{\infty}:= \sup\limits_{(x,y)\in \mathcal{X}^2} K(x,y)<+ \infty $,\,\,\,$\Vert  L\Vert_{\infty}:= \sup\limits_{(x,y)\in \mathcal{Y}^2} L(x,y)<+ \infty $;

\bigskip
\noindent it  ensures the existence of the kernel mean embeddings $m_X=\mathbb{E}\left(K(X,\cdot)\right)$,   $m_Y=\mathbb{E}\left( L(Y,\cdot)\right)$, and also that of $\mathbb{E}\left(L(Y,\cdot)\otimes K(X,\cdot)\right)$, where $\otimes$ denotes the tensor product defined as follows: for any $(f,g)\in\HX\times\HY$, $g\otimes f$ is the linear operator from $\HY$ to $\HX$ such that $(g\otimes f)(h)=\langle g,h\rangle_{\HY}\,f$, for any $h\in\HY$. 
For measuring dependence between $X$ and $Y$, Gretton et al. (2005) introduced the Hilbert-Schimdt Independence Criterion (HSIC)  defined as
\begin{equation}\label{hsic}
\textit{\textbf{H}}  = \bigg\Vert\mathbb{E} \big( L(Y,\cdot) \otimes K(X,\cdot) \big) - m_Y  \otimes m_X \bigg\Vert_{\textrm{HS}} ^2,\nonumber
\end{equation}
where $\Vert\cdot\Vert_{\textrm{HS}}$ denotes the Hilbert-Schmidt norm of operators. This permits to characterize independence and, consequently, to consider a test of independence between $X$ and $Y$, that is a test for the hypothesis $\mathscr{H}_0\,:\,\mathbb{P}_{XY}=\mathbb{P}_{X}\times\mathbb{P}_{Y}$ against  $\mathscr{H}_1\,:\,\mathbb{P}_{XY}\neq\mathbb{P}_{X}\times\mathbb{P}_{Y}$, where $\mathbb{P}_{XY}$ (resp. $\mathbb{P}_{X}$; resp. $\mathbb{P}_{Y}$) denotes the distribution of $(X,Y)$ (resp. $X$; resp. $Y$). Indeed, it is known from Theorem 4 of Gretton et al. (2005) that $\mathscr{H}_0$ is true  if and only if  $\textit{\textbf{H}}=0$ when the following assumption holds:

\bigskip

\noindent $(\mathscr{A}_2):$  $\mathcal{H}_X$ and  $\mathcal{H}_Y$ are  universal, and $\mathcal{X}$ and $\mathcal{Y}$  are  compact metric spaces.

\bigskip

\noindent Recall that a RKHS is universal if  it is dense, with respect to the $L_\infty$ norm,  in the space  of continuous maps    and if the associated kernel  is continuous (see, e.g., Gretton et al. (2012),  p. 727).  Therefore, a test of independence can be achieved by taking as test statistic a consistent estimator of  $\textit{\textbf{H}}$.  Letting $\{(X_i,Y_i)\}_{1\leq i\leq n}$ be a i.i.d. sample of $(X,Y)$ and replacing each expectation in (\ref{hsic}) by its empirical counterpart lead  to the estimator

\begin{eqnarray}\label{estemp}
\widehat{\textrm{\textbf{HSIC}}}_n&=& \bigg\Vert\frac{1}{n}\sum_{i=1}^n L(Y_i,\cdot) \otimes K(X_i,\cdot)   -\overline{L}_n  \otimes \overline{K}_n \bigg\Vert_{\textrm{HS}} ^2\nonumber\\
&= &\frac{1}{n^2}\sum_{i,j=1}^{n} k_{ij}\ell_{ij} + \frac{1}{n^4}\sum_{i,j,q,r=1}^{n} k_{ij}\ell_{qr} - \frac{2}{n^3} \sum_{i,j,q=1}^{n} k_{ij }\ell_{iq},
\end{eqnarray}
where $\overline{L}_n=n^{-1}\sum_{i=1}^n L(Y_i,\cdot)$,  $\overline{K}_n=n^{-1}\sum_{i=1}^n K(X_i,\cdot)$, 
$ k_{ij}=K(X_i,X_j)$ and  $\ell_{ij}=L(Y_i,Y_j)$.
This estimator was proposed  in Gretton (2005) as test statistic for testing for $\mathscr{H}_0$; but its  asymptotic distribution under null hypothesis    is an infinite sum of distributions   (see Zhang et al. (2018)). We propose an estimator obtained   from   an  appropriate modification of (\ref{estemp})   in order to yield asymptotic normality both under the  null hypothesis and under the alternative.
\section{Estimation of HSIC and asymptotic normality}
\noindent For $\gamma\in]0,1]$ and $r\in\mathbb{N}^\ast$, let $\left(w_{i,r}(\gamma)\right)_{1\leq i\leq r}$ be a sequence of  positive numbers satisfying:

\bigskip
\noindent $(\mathscr{A}_3):$ There exists  a strictly positive real number $ \tau $  and an  integer $ n_0 $ such that for all $r> n_0$:
$$r\left\vert\frac{1}{r}\sum_{i=1}^{r}w_{i,r}(\gamma)-1\right\vert\leq \tau.$$

\noindent $(\mathscr{A}_4):$ There exists  $C>0$ such that $\max\limits_{1\leq i \leq r}w_{i,r}(\gamma)<C$ for all  $r\in\mathbb{N}^\ast$ and $ \gamma \in]0, 1]$.

\noindent$(\mathscr{A}_5):$ For any   $ \gamma \in]0, 1]$,
$\lim\limits_{r\rightarrow +\infty}\frac{1}{r}\sum_{i=1}^{r}w^2_{i,r}(\gamma)=w^2(\gamma)> 1$.

 \bigskip

\noindent   A typical example is given by $w_{i,r}(\gamma)=1+(-1)^i\,\gamma$ (see Ahmad, 1993). We propose to estimate HSIC by a modification of (\ref{estemp}) given by

\begin{eqnarray}\label{stat}
\widehat{\textit{\textbf{H}}}_{n,\gamma}=\frac{1}{n^2}\sum_{i,j=1}^{n} k_{ij}\ell_{ij} + \frac{1}{n^4}\sum_{i,j,q,r=1}^{n} k_{ij}\ell_{qr} - \frac{2}{n^3} \sum_{i,j,q=1}^{n}w_{i,n}(\gamma)\, k_{ij }\ell_{iq}.\nonumber
\end{eqnarray}

\noindent Putting $
\mu=\mathbb{E}\left(  L(Y,\cdot)  \otimes K(X,\cdot)\right)$ and $ \nu=m_Y\otimes m_X$, 
and considering the functions  $\mathcal{U}$ and $\mathcal{V}$  from $\mathcal{X}\times\mathcal{Y}$ to $\mathbb{R}$ defined as
\begin{equation}\label{fonc1}
\begin{aligned}
 \mathcal{U} (x,y) &= \big\langle  L(y,\cdot)  \otimes K(x,\cdot) - \mu ,  \mu \big\rangle_{\textrm{HS}}\\
&+\big\langle  L(y,\cdot) \otimes m_X + m_Y \otimes  K(x,\cdot) - 2 \nu ,  \nu - \mu \big\rangle_{\textrm{HS}},\\
\mathcal{V}(x,y) &= \big\langle  L(y,\cdot)  \otimes K(x,\cdot)- \mu ,  \nu \big\rangle_{\textrm{HS}},
\end{aligned}
\end{equation}
where $\big\langle   \cdot , \cdot  \big\rangle_{\textrm{HS}}$  denotes the Hilbert-Smidt inner product,  we have:

\begin{thm}{}\label{thm}
Assume that $(\mathscr{A}_1)$ to  $(\mathscr{A}_5)$ hold. Then as $n\rightarrow +\infty$, one has  $ \sqrt{n}\Big (\widehat{\textit{\textbf{H}}}_{n,\gamma}-   \textit{\textbf{H}} \Big ) \stackrel{\mathscr D}{\longrightarrow} \mathcal N (0, \sigma ^2_\gamma)$, where 
\begin{equation*}
\begin{aligned}
\sigma ^2_\gamma   & =  4 Var\left(\mathcal{U}(X_1,Y_1)\right)+4 w^2(\gamma)Var\left(\mathcal{V}(X_1,Y_1)\right) -8Cov\left(\mathcal{U}(X_1,Y_1),\mathcal{V}(X_1,Y_1)\right).
\end{aligned}
\end{equation*}
\end{thm}
\bigskip

\noindent When  $\mathscr{H}_0$ is true, this theorem gives a   simpler expression for the variance of the limiting normal distribution. Indeed, one then has $\textit{\textbf{H}}=0$ and $\mu=\nu$. Therefore, $\mathcal{U}=\mathcal{V}$ and we have:

\begin{corl}{}\label{cor}
Assume that $(\mathscr{A}_1)$ to  $(\mathscr{A}_5)$ hold. Then, under  $\mathscr{H}_0$ ,  as $n\rightarrow +\infty$,   one has  $ \sqrt{n} \, \widehat{\textit{\textbf{H}}}_{n,\gamma} \stackrel{\mathscr D}{\longrightarrow} \mathcal N (0, \sigma ^2_\gamma)$, where 
$\sigma ^2_\gamma   =  4 \left(w(\gamma) ^2- 1\right)Var\left(\mathcal{V}(X_1,Y_1)\right)$.
\end{corl}
In order to achieve tests with the  above estimator, it is neessary to find a consistent estimator of the variance $\sigma ^2 _\lambda$. 
\begin{prop}\label{prop1}
Assume that $(\mathscr{A}_1)$ to  $(\mathscr{A}_5)$ hold. Then, under $\mathscr{H}_0$, the estimator $\widehat{\sigma}^2_\gamma=4 \left(w(\gamma) ^2- 1\right)\,\widehat{\alpha}$, where
\begin{equation}\label{sigma}
\begin{aligned}
\widehat{\alpha}&=
   \frac{1}{n}\sum_{i=1}^{n}  \Big ( \frac{1}{n}\sum_{j=1}^{n} \ell_{ij}k_{ij} - \frac{1}{n^2}\sum_{m=1}^{n} \sum_{p=1}^{n}\ell_{mp}k_{mp} \Big) ^2,\nonumber
\end{aligned}
\end{equation}
 is a consistent for $ \sigma ^2_\gamma$.
\end{prop}
\noindent The resulting test for independence is performed as follows: for a given significance level $\alpha\in]0,1[$, one has to reject $\mathscr{H}_0$ if  $\widehat{\textit{\textbf{H}}}_{n,\gamma}>n^{-1/2}\,\widehat{\sigma}_{\gamma}\,\Phi^{-1}(1-\alpha/2)$, where $\Phi$ is the cumulative distribution function of the standard normal distribution.
\section{Simulations}
\noindent In this section, we investigate the finite sample performance of the proposed test based on modified HSIC and compare it to  the test  of  Lai  et al. (2021)  based on angle covariance. For convenience, we denote our test as mhsic, and the test of   Lai  et al. (2021)  as   acov.   We computed   empirical sizes and powers    through   Monte Carlo simulations.  We considered  the case where $\mathcal{X}=\mathcal{Y}=\textrm{L}^2([0,1])$ and,  similar to the functional data considered in  Lai  et al. (2021), we take $X(t)=\sqrt{2}\sum_{k=1}^{50}\xi_k\,\cos(k\pi t)$ and  $Y(t)=\sqrt{2}\sum_{k=1}^{50}\nu_k\,\cos(k\pi t)$, where the $\xi_k$s are independent and distributed as the Cauchy distribution $\mathscr{C}(0,0.5)$ and, for a given $m\in\{0,\cdots,50\}$, $\nu_k=f(\xi_k)$ for $k=1,\cdots,m$ and the $\nu_k$s with $k=m+1,\cdots,50$ are sampled independently from the standard normal distribution. The null hypothesis $\mathscr{H}_0$ holds in case $m=0$, and the dependence level increases with $m$. Empirical sizes and powers were computed on the basis of $300$ independent  replicates. For each of them, we generated a sample  of size   $n=100$ of the above  processes in discretized versions  on equispaced values $t_1,\cdots,t_{51}$  in  $[0,1]$, where $t_j=(j-1)/50$, $j=1,\cdots,51$.  For performing  our method, we took $\gamma=0.32$ and used the gaussian  kernels $K(x,y)=L(x,y)=\exp\left(-\sigma\int_0^1\left( x(t)-y(t)\right)^2dt\right)$ with $\sigma^2=1/{150}$;   the terms $k_{ij}$ and $\ell_{ij}$
were computed by approximating   integrals involved in these   kernels  by using the trapezoidal rule. The significance level was taken as $\alpha=0.05$.  The  acov method  was used with $50$ permutations. Table 1 reports the obtained results. The obtained values for $m=0$  are close to the nominal size for all methods. For $m=1,3$, acov slightly outperforms our method which still remains competitive since the differences in the obtained values are low. For $m=5,10$, the two methods give high values for the power. This highlights the interest of the proposed test:  it is powerful enough and is fast compared to acov method which is based permutations and, therefore, leads   very high computation times.

\begin{table}
{\setlength{\tabcolsep}{0.08cm} 
\renewcommand{\arraystretch}{1} 
\begin{center}
{\begin{tabular}{cccccccccccccccccc}
\hline
 & &  &  & & &  &  &  &  &  &  &   &    &  & & \\
 & $f(x)$  & method &  &    & $m=0$ &  &  &$m=1$ &  &  &  $m=3$ & &   &  $m=5$&  & & $m=10$\\  

 &  &  & & &  &  & &  &   & &  &  &  &  &  & &  \\
\hline
\hline
 & & acov&   &   & 0.051&   &  & 1.00  &  &  &1.00  & &   &1.00  &   & &  1.00 \\
  & $x^3$& &   &    & &   &  &   &  &  &   & &   &   &   & &   \\
  & & mhsic&   &   & 0.052&   &  & 0.89  &  &  &0.92  & &   &1.00  &   & &  1.00 \\
\hline
 & & acov&   &   & 0.053&   &  & 0.84  &  &  &0.97  & &   &1.00  &   & &  1.00 \\
  & $x^2$& &   &    & &   &  &   &  &  &   & &   &   &   & &   \\
  & & mhsic&   &   & 0.060&   &  & 0.80  &  &  &0.91  & &   &0.94  &   & &  1.00 \\
\hline
 & & acov&   &   & 0.051&   &  & 0.85  &  &  &0.95  & &   &0.95  &   & &  0.97 \\
  & $x^2\sin(x)$& &   &    & &   &  &   &  &  &   & &   &   &   & &   \\
  & & mhsic&   &   & 0.051&   &  & 0.80  &  &  &0.88  & &   &0.96  &   & &  0.99 \\
\hline\hline
\end{tabular}}         
\end{center}}
\centering \caption{\label{table:tab4}Empirical sizes and powers over 300 replications  with significance level $\alpha=0.05$.}
\end{table}
\section{Proofs}
\subsection{A technical lemma}
\begin{lem}{}
Let $(Z_i)_{1 \leq i \leq n}$ be an i.i.d.  sample of a random variable  $Z$ valued into  a Hilbert space $\mathcal{H}$ and such that $\mathbb{E}\left(\Vert Z\Vert_\mathcal{H}^2\right)<+\infty$. Then
\begin{equation}\label{eqt2}
\begin{aligned}
\Big\Vert\frac{1}{n} \sum_{i=1}^{n} \left(w_{i,n}(\gamma) - 1 \right) Z_i \Big\Vert_{\mathcal{H}} = o_p (1).
\end{aligned}
\end{equation}
\end{lem}

\noindent\textit{ Proof.} Putting $m_Z=\mathbb{E}(Z)$, we have:
\begin{equation*}
\begin{aligned}
& \Big\Vert  \frac{1}{n} \sum_{i=1}^{n}\left(w_{i,n}(\gamma) - 1 \right) Z_i  \Big\Vert  _{\mathcal{H}}\\
& = \Big\Vert  \frac{1}{n} \sum_{i=1}^{n}w_{i,n}(\gamma)\left( Z_i - m_Z\right)  - \frac{1}{n} \sum_{i=1}^{n} Z_i + m_Z +  \frac{1}{n} \sum_{i=1}^{n}\left(w_{i,n}(\gamma) - 1 \right) m_Z\Big\Vert  _{\mathcal{H}} \\
& \leq \Big\Vert  \frac{1}{n} \sum_{i=1}^{n}w_{i,n}(\gamma)\left( Z_i - m_Z\right) \Big\Vert  _{\mathcal{H}} + \Big\Vert  \frac{1}{n} \sum_{i=1}Z_i - m_Z \Big\Vert  _{\mathcal{H}} + \Big|\frac{1}{n}\sum_{i=1}^{n}\left(w_{i,n}(\gamma) - 1 \right)\Big| \times \Big\Vert  m_Z\Big\Vert  _{\mathcal{H}}.
\end{aligned}
\end{equation*}
First, using the equality
\begin{equation*}
\begin{aligned}
\mathbb{E}\left(  \Big\Vert  \frac{1}{n} \sum_{i=1}^{n}w_{i,n}(\gamma)\left(Z_i - m_Z\right) \Big\Vert  _{\mathcal{H}}  ^2 \right)= \left(\frac{1}{n^2} \sum_{i=1}^{n}w_{i,n}^2 (\gamma)\right) \mathbb{E}\left( \Big\Vert  Z - m_Z \Big\Vert  _{\mathcal{H}}  ^2\right),
\end{aligned}
\end{equation*}
assumption $(\mathscr{A}_5)$ and Markov inequality, we get $ \Big\Vert  \frac{1}{n} \sum_{i=1}^{n}w_{i,n}(\gamma)\left( Z_i - m_Z\right) \Big\Vert  _{\mathcal{H}}=o_p(1)$.  Secondly,   the law of large numbers gives $ \Big\Vert  \frac{1}{n} \sum_{i=1}Z_i - m_Z \Big\Vert  _{\mathcal{H}} =o_p(1)$. Thirdly, using  assumption $(\mathscr{A}_3)$ we obtain  $ \Big|\frac{1}{n}\sum_{i=1}^{n}\left(w_{i,n}(\gamma) - 1 \right)\Big| \rightarrow 0$ as $n\rightarrow +\infty$. Then,  (\ref{eqt2}) is obtained.

\subsection{Proof of Theorem 1}
\noindent Using the definition of the Hilbert-Schmidt inner product, the reproducing properties of $K$ and $L$, and the equality $(a\otimes b)^\ast=b\otimes a$, it is easy to see that
\begin{eqnarray}
\widehat{\textit{\textbf{H}}}_{n,\gamma}&=&\bigg\Vert\frac{1}{n}\sum_{i=1}^n L(Y_i,\cdot) \otimes K(X_i,\cdot) \bigg\Vert_{\textrm{HS}} ^2 +\bigg\Vert\overline{L}_n  \otimes \overline{K}_n \bigg\Vert_{\textrm{HS}} ^2\nonumber\\
& & - \frac{2}{n} \sum_{i=1}^{n}w_{i,n}(\gamma)\Big\langle L(Y_i,\cdot) \otimes K(X_i,\cdot),  \overline{L}_n  \otimes \overline{K}_n \Big\rangle_{\textrm{HS}}\nonumber
\end{eqnarray}
and
\begin{eqnarray}
\textit{\textbf{H}}&=&\bigg\Vert \mathbb{E}\Big( L(Y,\cdot) \otimes K(X,\cdot)\Big) \bigg\Vert_{\textrm{HS}} ^2 +\bigg\Vert m_Y  \otimes m_X \bigg\Vert_{\textrm{HS}} ^2-2\Big\langle\mathbb{E}\Big( L(Y,\cdot) \otimes K(X,\cdot)\Big),   m_Y  \otimes m_X \Big\rangle_{\textrm{HS}}\nonumber.
\end{eqnarray}
Therefore, using the two following equalities: $\Vert a\Vert^2=\Vert a-b\Vert^2+2\langle a,b\rangle-\Vert b\Vert^2$, $\langle a_1,a_2\rangle=\langle a_1-b_1,a_2-b_2\rangle+\langle a_1,b_2\rangle+\langle a_2,b_1\rangle-\langle b_1,b_2\rangle$, and putting $k_{X_i}=K(X_i,\cdot)$, 
  $l_{Y_i}=L(Y_i,\cdot)$, 
$ \overline{l \otimes k}^n =\frac{1}{n} \sum_{i=1}^{n}  l_{Y_i} \otimes k_{X_i} $,  we get
\begin{equation}\label{decompo}
\begin{aligned}
&   \sqrt{n}\Big (\widehat{\textit{\textbf{H}}}_{n,\gamma} - \textit{\textbf{H}}  \Big )\nonumber \\
& = \sqrt{n}\Big [\,\, \Big\Vert  \overline{l \otimes k}^n - \mu\Big\Vert  _{\textrm{HS}}  ^2 + 2 \; \big\langle \overline{l \otimes k}^n,  \mu  \big\rangle_{\textrm{HS}} - \Big\Vert   \mu \Big\Vert  _{\textrm{HS}}  ^2  + \Big\Vert  \overline{L}_n \otimes \overline{K}_n - \nu\Big\Vert  _{\textrm{HS}}  ^2 + 2 \; \big\langle \overline{L}_n \otimes \overline{K}_n ,  \nu  \big\rangle_{\textrm{HS}} - \Big\Vert   \nu \Big\Vert  _{\textrm{HS}}  ^2\nonumber \\
&  \qquad   - \frac{2}{n} \sum_{i=1}^{n} \left(w_{i,n}(\gamma) - 1 \right) \;  \big\langle  l_{Y_i} \otimes k_{X_i} , \overline{L}_n \otimes \overline{K}_n  - \nu \big\rangle_{\textrm{HS}} - \frac{2}{n} \sum_{i=1}^{n}w_{i,n}(\gamma) \;  \big\langle  l_{Y_i} \otimes k_{X_i} , \nu \big\rangle_{\textrm{HS}}\nonumber\\
&  \qquad  +  \frac{2}{n} \sum_{i=1}^{n}\;  \big\langle  l_{Y_i} \otimes k_{X_i} ,  \nu \big\rangle_{\textrm{HS}} - 2  \big\langle \overline{l \otimes k}^n - \mu , \overline{L}_n \otimes \overline{K}_n  - \nu \big\rangle_{\textrm{HS}} - \frac{2}{n} \sum_{i=1}^{n} \;  \big\langle  l_{Y_i} \otimes k_{X_i} ,  \nu \big\rangle_{\textrm{HS}}  \nonumber \\
&   \qquad   - 2 \; \big\langle \overline{L}_n \otimes \overline{K}_n ,  \mu  \big\rangle_{\textrm{HS}}- \frac{2}{n} \sum_{i=1}^{n} \left(w_{i,n}(\gamma) - 1 \right) \;  \big\langle \mu ,  \nu \big\rangle_{\textrm{HS}} +\frac{2}{n} \sum_{i=1}^{n}w_{i,n}(\gamma) \;  \big\langle \mu ,  \nu \big\rangle_{\textrm{HS}} - \Big\Vert   \mu - \nu \Big\Vert  _{\textrm{HS}}  ^2\,\,\, \Big]\nonumber\\
&= A_n+B_n+C_n+D_n,
\end{aligned}
\end{equation}
where
\begin{equation}
\begin{aligned}
A_n &=\sqrt{n}\Big( \Big\Vert  \overline{l \otimes k}^n - \mu\Big\Vert  _{\textrm{HS}}  ^2+ \Big\Vert  \overline{L}_n \otimes \overline{K}_n - \nu\Big\Vert  _{\textrm{HS}}  ^2\Big)\\
& = n^{-1/2} \bigg[\,\,\Big\Vert  \sqrt{n}\big(\overline{l \otimes k}^n - \mu \big)\Big\Vert  _{\textrm{HS}}  ^2+ \Big\Vert  \left(\sqrt{n}(\overline{L}_n - m_Y)\right)\otimes\left(\overline{K}_n - m_X\right)\nonumber \\
& +\left(\sqrt{n}(\overline{L}_n - m_Y)\right)\otimes m_X + m_Y \otimes\left(\sqrt{n}(\overline{K}_n - m_X)\right) \Big\Vert  _{\textrm{HS}}  ^2\,\,\bigg],\\
\end{aligned}
\end{equation}
\begin{equation}
\begin{aligned}
B_n &= -2 \bigg [     \Big\langle \frac{1}{n} \sum_{i=1}^{n} \left(w_{i,n}(\gamma) - 1 \right) l_{Y_i} \otimes k_{X_i} , \sqrt{n} \big(\overline{L}_n \otimes \overline{K}_n  - \nu \big) \Big\rangle_{\textrm{HS}} + \frac{\sqrt{n}}{n} \sum_{i=1}^{n} \left(w_{i,n}(\gamma) - 1 \right) \;  \big\langle \mu ,  \nu \big\rangle_{\textrm{HS}}\bigg],\nonumber\\
C_n &=- 2  \Big\langle \sqrt{n} \big(\overline{l \otimes k}^n - \mu \big), \overline{L}_n \otimes \overline{K}_n  - \nu \Big\rangle_{\textrm{HS}}+ 2\big\langle \sqrt{n}(\overline{L}_n - m_Y) \otimes (\overline{K}_n - m_X),  \nu - \mu \big\rangle_{\textrm{HS}},\nonumber \\
D_n &=  \frac{2}{\sqrt{n}}\sum_{i=1}^{n}\bigg( \mathcal{U} (X_i,Y_i) - w_{i,n}(\gamma)\mathcal{V}  (X_i,Y_i) \bigg), 
\end{aligned}
\end{equation}
$\mathcal{U}$ and $\mathcal{V}$ being the functions defined in (\ref{fonc1}). From the central limit theorem,  $ \sqrt{n}\big(\overline{l \otimes k}^n -\mu)$, $\sqrt{n}(\overline{L}_n - m_Y)$ and $\sqrt{n}(\overline{K}_n - m_X)$ converge in distribution, as $n\rightarrow +\infty$, to random variables having normal distributions. In addition, by the law of large numbers, $\overline{K}_n - m_X$  converges  in probability to $0$, as $n\rightarrow +\infty$. We then deduce that $A_n=o_p(1)$. On the other hand, using Cauchy-Schwartz inequality, we obtain
\begin{eqnarray}
\vert B_n\vert &\leq& 2 \bigg [    \Big\Vert  \frac{1}{n} \sum_{i=1}^{n} \left(w_{i,n}(\gamma) - 1 \right) l_{Y_i} \otimes k_{X_i} \Big\Vert_{\textrm{HS}}\Big\Vert \sqrt{n} \big(\overline{L}_n \otimes \overline{K}_n  - \nu \big) \Big\Vert_{\textrm{HS}} \nonumber\\
& &+\sqrt{n}\,\Big\vert \frac{1}{n} \sum_{i=1}^{n}  w_{i,n}(\gamma) - 1 \Big\vert \,\,\,  \big\Vert\mu \big\Vert_{\textrm{HS}}\,\, \big\Vert\  \nu  \big\Vert_{\textrm{HS}}\bigg].\nonumber
\end{eqnarray}
Snce $\sqrt{n} \big(\overline{L}_n \otimes \overline{K}_n  - \nu\big)=\big(\sqrt{n} \big(\overline{L}_n -m_Y\big)\big)\otimes \overline{K}_n  +m_Y \otimes\big(\sqrt{n} \big( \overline{K}_n  - m_X\big)\big)$, and from  $(\mathscr{A}_3)$, $\Big\vert \frac{1}{n} \sum_{i=1}^{n}  w_{i,n}(\gamma) - 1 \Big\vert \leq \tau/n$ for $n$ large enough,  it follows:
\begin{eqnarray}
\vert B_n\vert &\leq& 2 \bigg [    \Big\Vert  \frac{1}{n} \sum_{i=1}^{n} \left(w_{i,n}(\gamma) - 1 \right) l_{Y_i} \otimes k_{X_i} \Big\Vert_{\textrm{HS}}\bigg(\Big\Vert\big(\sqrt{n} \big(\overline{L}_n -m_Y\big)\big)\otimes \overline{K}_n  \big) \Big\Vert_{\textrm{HS}} \nonumber\\
& &+\Big\Vert m_Y \otimes\big(\sqrt{n} \big( \overline{K}_n  - m_X\big)\big)\Big\Vert_{\textrm{HS}}\bigg) +\frac{\tau}{\sqrt{n}} \,\,\,  \big\Vert\mu \big\Vert_{\textrm{HS}}\,\, \big\Vert\  \nu  \big\Vert_{\textrm{HS}}\bigg].\nonumber\\
&=& 2 \bigg [    \Big\Vert  \frac{1}{n} \sum_{i=1}^{n} \left(w_{i,n}(\gamma) - 1 \right) l_{Y_i} \otimes k_{X_i} \Big\Vert_{\textrm{HS}}\bigg(\Big\Vert\sqrt{n} \big(\overline{L}_n -m_Y\big)\Big\Vert_{\mathcal{H}_Y}\,\Big\Vert \overline{K}_n   \Big\Vert_{\mathcal{H}_X} \nonumber\\
& &+\Big\Vert m_Y \Big\Vert_{\mathcal{H}_Y}\Big\Vert\sqrt{n} \big( \overline{K}_n  - m_X\big)\Big\Vert_{\mathcal{H}_X}\bigg) +\frac{\tau}{\sqrt{n}} \,\,\,  \big\Vert\mu \big\Vert_{\textrm{HS}}\,\, \big\Vert\  \nu  \big\Vert_{\textrm{HS}}\bigg].\nonumber
\end{eqnarray}
Further, using the reproducing property  of $K$, we obtain
\begin{equation}\label{ineg}
\Big\Vert \overline{K}_n   \Big\Vert_{\mathcal{H}_X}\leq \frac{1}{n}\sum_{i=1}^n\Big\Vert K(X_i,\cdot)   \Big\Vert_{\mathcal{H}_X}= \frac{1}{n}\sum_{i=1}^n\sqrt{ K(X_i,X_i) }\leq\Vert K\Vert_\infty^{1/2}. 
\end{equation}
Hence
\begin{eqnarray}
\vert B_n\vert &\leq& 2 \bigg [    \Big\Vert  \frac{1}{n} \sum_{i=1}^{n} \left(w_{i,n}(\gamma) - 1 \right) l_{Y_i} \otimes k_{X_i} \Big\Vert_{\textrm{HS}}\bigg(\Big\Vert\sqrt{n} \big(\overline{L}_n -m_Y\big)\Big\Vert_{\mathcal{H}_Y}\,\Vert K\Vert_\infty^{1/2} \nonumber\\
& &+\Big\Vert m_Y \Big\Vert_{\mathcal{H}_Y}\Big\Vert\sqrt{n} \big( \overline{K}_n  - m_X\big)\Big\Vert_{\mathcal{H}_X}\bigg) +\frac{\tau}{\sqrt{n}} \,\,\,  \big\Vert\mu \big\Vert_{\textrm{HS}}\,\, \big\Vert\  \nu  \big\Vert_{\textrm{HS}}\bigg].\nonumber
\end{eqnarray}
From Lemma 1, we have $ \Big\Vert  \frac{1}{n} \sum_{i=1}^{n} \left(w_{i,n}(\gamma) - 1 \right) l_{Y_i} \otimes k_{X_i} \Big\Vert_{\textrm{HS}}=o_p(1)$. Then, since $\sqrt{n} \big(\overline{K}_n -m_X\big)$ and $\sqrt{n} \big(\overline{L}_n -m_Y\big)$ converge in distribution as $n\rightarrow +\infty$, we deduce from the preceding inequality that $B_n=o_p(1)$. Using again Cauchy-Scwartz inequality,  we get
\begin{equation}
\begin{aligned}
\vert C_n \vert&\leq  2   \Big\Vert  \sqrt{n} \big(\overline{l \otimes k}^n - \mu  \big)\Big\Vert _{\textrm{HS}}   \Big\Vert  \overline{L}_n \otimes \overline{K}_n  - \nu  \Big\Vert _{\textrm{HS}}+ 2\Big\Vert \sqrt{n}(\overline{L}_n - m_Y) \otimes (\overline{K}_n - m_X)\Big\Vert _{\textrm{HS}}\Big\Vert   \nu - \mu \Big\Vert _{\textrm{HS}}\nonumber \\
&= 2   \Big\Vert  \sqrt{n} \big(\overline{l \otimes k}^n - \mu  \big)\Big\Vert _{\textrm{HS}}   \Big\Vert  \overline{L}_n \otimes \overline{K}_n  - \nu  \Big\Vert _{\textrm{HS}}+ 2\Big\Vert \sqrt{n}(\overline{L}_n - m_Y) \Big\Vert_{\mathcal{H}_Y}\Big\Vert \overline{K}_n - m_X\Big\Vert_{\mathcal{H}_X} \Big\Vert   \nu - \mu \Big\Vert _{\textrm{HS}}.\nonumber
\end{aligned}
\end{equation}
We already know that, as $n\rightarrow +\infty$,   $\overline{l \otimes k}^n$ and $\sqrt{n}(\overline{L}_n - m_Y)$ converge in distribution to normal random variables,  and that $\overline{K}_n$ and $\overline{L}_n$ converge in probability to  $m_X$ and $m_Y$ respectively. Thus $ \overline{L}_n \otimes \overline{K}_n$ converges in probability to $\nu$ as $n\rightarrow +\infty$, and the preceding inequality implies that $C_n=o_p(1)$. We can conclude that
\begin{equation}\label{eqt5}
\begin{aligned}
 \sqrt{n}\Big (\widehat{\textit{\textbf{H}}}_{n,\gamma} - \textit{\textbf{H}}  \Big ) = D_n+o_p(1)= \frac{2}{\sqrt{n}}\sum_{i=1}^{n}\bigg( \mathcal{U} (X_i,Y_i) - w_{i,n}(\gamma)\mathcal{V}  (X_i,Y_i) \bigg)  + o_p (1)  \nonumber
\end{aligned}
\end{equation}
and, consequently, that $\sqrt{n}\Big (\widehat{\textit{\textbf{H}}}_{n,\gamma} - \textit{\textbf{H}}  \Big ) $ has the same limiting distribution than $D_n$; it remains to derive this latter.  Let us set  
\[
s_{n,\gamma}^2=\sum_{i=1}^{n}Var\Big(\mathcal{U}(X_i,Y_i)-w_{i,n}(\gamma)\mathcal{V}(X_i,Y_i)\Big).
\] 
 By similar arguments as in the proof of Theorem 1 in Makigusa and Naito (2020) we obtain that, for any $\varepsilon>0$,
\[
s_{n,\gamma}^{-2}\sum_{i=1}^{n}\int_{\{(x,y):| \mathcal{U}(x,y)-w_{i,n}(\gamma)\mathcal{V}(x,y) |>\varepsilon s_{n,\gamma}\}}^{} \bigg(\mathcal{U}(x,y)-w_{i,n}(\gamma)\mathcal{V}(x,y)\bigg)^2\,d\mathbb{P}_{XY}(x,y) 
\]
converges to $0$ as $n\rightarrow +\infty$. Therefore, by Section 1.9.3 in Serfling (1980) we obtain that  
$\sqrt{n}s_{n,\gamma}^{-1}\frac{D_n}{2}	\stackrel{\mathscr{D}}{\rightarrow} \mathcal{N}\left(0,1\right)$. However,
\begin{eqnarray*}
\left(\frac{s_{n,\gamma}}{\sqrt{n}}\right)^2
&=&Var\left(\mathcal{U}(X_1,Y_1)\right)+\left(\frac{1}{n}\sum_{i=1}^{n}w^2_{i,n}(\gamma)\right)Var\left(\mathcal{V}(X_1,Y_1)\right)\nonumber\\
&&-2\left(\frac{1}{n}\sum_{i=1}^{n}w_{i,n}(\gamma)\right)Cov\left(\mathcal{U}(X_1,Y_1),\mathcal{V}(X_1,Y_1)\right),
\end{eqnarray*}
from  $ (\mathscr{A}_3) $ and $ (\mathscr{A}_5) $, 
\[
\lim_{n\rightarrow +\infty}\left(\frac{1}{n}\sum_{i=1}^{n}w^2_{i,n}(\gamma)\right)=w^2(\gamma)\,\,\,\textrm{ and }\,\,\,\lim_{n\rightarrow +\infty}\left(\frac{1}{n}\sum_{i=1}^{n}w_{i,n}(\gamma)\right)=1.
\]
Thus
\[
\lim\limits_{n\rightarrow +\infty}\left(n^{-1}s_{n,\gamma}^2\right)=Var\left(\mathcal{U}(X_1,Y_1)\right)+ w^2(\gamma)Var\left(\mathcal{V}(X_1,Y_1)\right) -2Cov\left(\mathcal{U}(X_1,Y_1),\mathcal{V}(X_1,Y_1)\right)
\]
and, therefore, $D_n	\stackrel{\mathscr{D}}{\rightarrow} \mathcal{N}\left(0,\sigma_\gamma^2\right)$. 

\subsection{Proof of Proposition $1$}
\noindent From  the definition of $\big\langle\cdot , \cdot\big\rangle_{\textrm{HS}}$ and the reproducing properties of $K$ and $L$ it is easilly seen that  $ \big\langle  l_{Y_i} \otimes k_{X_i} ; \overline{l \otimes k}^n \big\rangle_{\textrm{HS}}=n^{-1}\sum_{j=1}^{n} \ell_{ij}k_{ij}$. Hence
\begin{eqnarray}
\widehat{\alpha}&=&
   \frac{1}{n}\sum_{i=1}^{n}  \bigg ( \big\langle  l_{Y_i} \otimes k_{X_i} , \overline{l \otimes k}^n \big\rangle_{\textrm{HS}}- \frac{1}{n}\sum_{m=1}^{n} \big\langle  l_{Y_m} \otimes k_{X_m} , \overline{l \otimes k}^n \big\rangle_{\textrm{HS}} \bigg) ^2\nonumber\\
&=&
   \frac{1}{n}\sum_{i=1}^{n}  \big\langle  l_{Y_i} \otimes k_{X_i} , \overline{l \otimes k}^n \big\rangle_{\textrm{HS}}^2-  \bigg (\frac{1}{n}\sum_{i=1}^{n} \big\langle  l_{Y_i} \otimes k_{X_i} , \overline{l \otimes k}^n \big\rangle_{\textrm{HS}} \bigg) ^2\nonumber.
\end{eqnarray}
Let us notice that
\begin{equation}\label{toUse3}
\begin{aligned}
& \frac{1}{n}\sum_{i=1}^{n}   \big\langle  l_{Y_i} \otimes k_{X_i} , \overline{l \otimes k}^n \big\rangle_{\textrm{HS}}^2 - \frac{1}{n}\sum_{i=1}^{n}   \big\langle  l_{Y_i} \otimes k_{X_i} , \mu \big\rangle_{\textrm{HS}}^2\\
& =  \frac{1}{n}\sum_{i=1}^{n}   \big\langle  l_{Y_i} \otimes k_{X_i} , \overline{l \otimes k}^n - \mu \big\rangle_{\textrm{HS}}^2 + \frac{2}{n}\sum_{i=1}^{n} \big\langle  l_{Y_i} \otimes k_{X_i} , \mu \big\rangle_{\textrm{HS}}  \big\langle  l_{Y_i} \otimes k_{X_i} , \overline{l \otimes k}^n - \mu \big\rangle_{\textrm{HS}}\nonumber.
\end{aligned}
\end{equation}
Then  Cauchy-Schwartz inequality,  reproducing property  and assumption $(\mathscr{A}_1)$ give
$$ \Big|  \frac{1}{n}\sum_{i=1}^{n}   \big\langle  l_{Y_i} \otimes k_{X_i} , \overline{l \otimes k}^n - \mu \big\rangle_{\textrm{HS}}^2 \Big | \leq\Vert K\Vert_{\infty} \Vert L\Vert_{\infty}\,\Big\Vert   \overline{l \otimes k}^n - \mu  \Big\Vert_{\textrm{HS}}  ^2   $$
and
$$    \Big| \frac{1}{n}\sum_{i=1}^{n} \big\langle  l_{Y_i} \otimes k_{X_i} , \mu \big\rangle_{\textrm{HS}}  \big\langle  l_{Y_i} \otimes k_{X_i} , \overline{l \otimes k}^n - \mu \big\rangle_{\textrm{HS}} \Big | \leq  \Big\Vert K \Big\Vert_{\infty} \Big\Vert L \Big\Vert_{\infty} \Big\Vert\mu \Big\Vert_{\textrm{HS}} \Big\Vert   \overline{l \otimes k}^n - \mu  \Big\Vert_{\textrm{HS}}.   $$
Since  $ \Big\Vert   \overline{l \otimes k}^n - \mu  \Big\Vert _{\textrm{HS}}   = o_p(1) $,  these inequalities show that 
$$ \frac{1}{n}\sum_{i=1}^{n}   \big\langle  l_{Y_i} \otimes k_{X_i} , \overline{l \otimes k}^n \big\rangle_{\textrm{HS}}^2 - \frac{1}{n}\sum_{i=1}^{n}   \big\langle  l_{Y_i} \otimes k_{X_i} , \mu \big\rangle_{\textrm{HS}}^2 = o_p(1) .$$
The law of large numbers and Slutsky's theorem allow to  conclude that the sequence  $\frac{1}{n}\sum_{i=1}^{n}   \big\langle  l_{Y_i} \otimes k_{X_i},  \overline{l \otimes k}^n \big\rangle_{\textrm{HS}}^2 $ converges in probabilty to $\mathbb{E} \left(\big\langle  l_{Y} \otimes k_{X} , \mu \big\rangle_{\textrm{HS}}^2\right)$ as $n\rightarrow +\infty$.
Similarly, from
\begin{equation}\label{toUse2}
\begin{aligned}
& \Big\vert\frac{1}{n}\sum_{i=1}^{n}   \big\langle  l_{Y_i} \otimes k_{X_i} , \overline{l \otimes k}^n \big\rangle_{\textrm{HS}}    - \frac{1}{n}\sum_{i=1}^{n}   \big\langle  l_{Y_i} \otimes k_{X_i} , \mu \big\rangle_{\textrm{HS}} \Big\vert \nonumber  \\
& =\Big\vert  \frac{1}{n}\sum_{i=1}^{n}   \big\langle  l_{Y_i} \otimes k_{X_i} , \overline{l \otimes k}^n - \mu \big\rangle_{\textrm{HS}}\Big\vert \leq  \Big\Vert K\Big\Vert_{\infty}^{1/2} \Big\Vert L\Big\Vert_{\infty}^{1/2}\Big\Vert   \overline{l \otimes k}^n - \mu  \Big\Vert  _{\textrm{HS}}\nonumber
\end{aligned}
\end{equation}
we get  $n^{-1}\sum_{i=1}^{n}   \big\langle  l_{Y_i} \otimes k_{X_i} , \overline{l \otimes k}^n \big\rangle_{\textrm{HS}}    - n^{-1}\sum_{i=1}^{n}   \big\langle  l_{Y_i} \otimes k_{X_i} , \mu \big\rangle_{\textrm{HS}} =o_p(1)$,  and   $n^{-1}\sum_{i=1}^{n}   \big\langle  l_{Y_i} \otimes k_{X_i} , \overline{l \otimes k}^n \big\rangle_{\textrm{HS}}$ converges in probability to $\mathbb{E}\left( \big\langle  l_{Y} \otimes k_{X} ,  \mu \big\rangle_{\textrm{HS}}\right)$. So,  $\widehat{\alpha}$ is a consistent estimator of $Var\left( \big\langle  l_{Y_1} \otimes k_{X_1},  \mu \big\rangle_{\textrm{HS}}\right)=Var\left(\mathcal{V}(X_1,Y_1)\right)$.

\end{document}